\documentclass[11pt]{amsart}
\usepackage{amsmath,amssymb,txfonts}
\usepackage{amssymb}
\usepackage{amsmath}
\usepackage{mathrsfs}
\usepackage[shortlabels]{enumitem}
\usepackage{color}
\usepackage{extarrows}
\newtheorem{theorem}{Theorem}[section]
\newtheorem{lemma}[theorem]{Lemma}
\newtheorem{proposition}[theorem]{Proposition}

\theoremstyle{definition}
\newtheorem{definition}[theorem]{Definition}

\theoremstyle{remark}

\numberwithin{equation}{section}

\begin{document}

\title[Fractional Besov Spaces]{Fractional Besov Trace/Extension Type Inequalities via the Caffarelli-Silvestre extension}

\author[P. Li]{Pengtao Li}
\address[P. Li]{College of Mathematics, Qingdao University, Qingdao, Shandong 266071, China}
\email{ptli@qdu.edu.cn}

\author[R. Hu]{Rui Hu}
\address[R. Hu]{Department of Mathematics and Statistics,  MacEwan University   Edmonton, Alberta T5J2P2 Canada}
\email{hur3@macewan.ca}

\author[Z. Zhai]{Zhichun Zhai}
\address[Z. Zhai]{Department of Mathematics and Statistics,  MacEwan University,   Edmonton, Alberta T5J2P2 Canada
}
\email{zhaiz2@macewan.ca}

\thanks{Corresponding author:  zhaiz2@macewan.ca}
\thanks{Project supported:  Pengtao Li was  supported by National Natural Science Foundation of China (No. 11871293, No. 12071272) and Shandong Natural Science Foundation of China (No. ZR2017JL008).}

\date{}

\dedicatory{}


\begin{abstract} Let $u(\cdot,\cdot)$ be the Caffarelli-Silvestre extension of $f.$
The first goal of  this article is to  establish the fractional  trace type  inequalities involving the Caffarelli-Silvestre extension $u(\cdot,\cdot)$ of $f.$ In doing so, firstly, we  establish the
fractional Sobolev/ logarithmic Sobolev/ Hardy trace inequalities in terms of $\nabla_{(x,t)}u(x,t).$ Then, we prove the  fractional anisotropic Sobolev/ logarithmic Sobolev/ Hardy trace inequalities  in terms of $ {\partial_{t} u(x,t)}$ or  $(-\Delta)^{-\gamma/2}u(x,t)$ only.
Moreover, based on   an estimate of the Fourier transform of the Caffarelli-Silvestre extension kernel and the sharp affine weighted $L^p$ Sobolev inequality, we prove that the $\dot{H}^{-\beta/2}(\mathbb{R}^n)$ norm of $f(x)$ can be controlled by the product of the weighted $L^p-$affine energy and the  weighted $L^p-$norm of ${\partial_{t} u(x,t)}.$
The second goal of this article is to
 characterize non-negative  measures $\mu$ on $\mathbb{R}^{n+1}_+$ such that the embeddings $$\|u(\cdot,\cdot)\|_{L^{q_0,p_0}_{\mu}(\mathbb{R}^{n+1})}\lesssim \|f\|_{\dot{\Lambda}^{p,q}_\beta(\mathbb{R}^n)}$$
 hold for some $p_0$ and $q_0$ depending on $p$ and $q$  which are classified in  three different cases: (1).    $p=q\in (n/(n+\beta),1];$ (2) $(p,q)\in (1,n/\beta)\times (1,\infty);$  (3).  $(p,q)\in (1,n/\beta)\times\{\infty\}.$
For case (1),
the embeddings can be characterized in terms of an analytic condition of the variational capacity minimizing function,  the iso-capacitary inequality of open balls, and other weak type inequalities.
For cases (2) and (3), the embeddings  are characterized  by the iso-capacitary inequality for fractional Besov capacity of open sets.
\end{abstract}

\maketitle

\tableofcontents \pagenumbering{arabic}


\section{Introduction}
Sobolev type  inequalities  are extremely important in analysis and geometry,  and have been applied in studying of partial differential equations. In this article, we will establish factional Besov trace and extension type inequalities involving fractional Laplacian and the Caffarelli-Silvestre extension.

For $s\in (0,2),$ the fractional Laplace operator $(-\Delta)^{s/2} $
in $\mathbb{R}^n $ is defined  on the Schwartz class through the Fourier transform as

$$\widehat{[(-\Delta)^{s/2} f]} (\xi)= |\xi|^{s}\widehat{f}(\xi),$$
or via the Riesz potential as
$$
(-\Delta)^{s/2} f(x)=\frac{s2^{s}\Gamma({(n+s)}/{2})}{2\Gamma(1-s/2)\pi^{{n}/{2}}}\hbox{P.V.}\int_{\mathbb{R}^n}\frac{f(x)-f(y)}{|x-y|^{n+s}}dy,$$
where $\widehat{f}(\xi)=\int_{\mathbb{R}^n}e^{- ix\cdot \xi } f(x) dx$ is the Fourier transform of $f,$ $\hbox{P.V.}$ denotes the Cauchy principle value integral.
The case $s=2$ becomes the standard local Laplacian.   The general case can be written  in a local way by the so-called Caffarelli-Silvestre extension introduced by  Caffarelli and  Silvestre  in \cite{Caffarelli}.

Let $f$ be a regular function in $\mathbb{R}^n.$ We say that $u=u(\cdot,\cdot)$ is the Caffarelli-Silvestre extension of $f$ to the upper half-space $\mathbb{R}^{n+1}_{+}:=\mathbb{R}^n\times (0,\infty),$ if $u$ is a solution to the problem
\begin{equation}\label{eq1}
\left\{
\begin{aligned}
\hbox{div}(t^{1-s}\nabla u)&= 0,\quad \hbox{in}\quad \mathbb{R}^{n+1}_{+};\\
u &= f, \quad \hbox{on} \quad \mathbb{R}^n\times\{t=0\}.
\end{aligned}
\right.
\end{equation}
The Caffarelli-Silvestre extension is well defined for smooth functions through the  fractional Poisson  kernel: $$p^{s}_t(x):=\frac{c(n,s)t^s}{(|x|^2+t^2)^{{(n+s)}/{2}}},\quad\quad c(n,s)=\frac{\Gamma({(n+s)}/{2})}{\pi^{n/2}\Gamma({s}/{2})}$$
 as follows
$$u(x,t)=p_t^{s}\ast f(x)=c(n,s)t^s\int_{\mathbb{R}^n}\frac{f(y)}{(|x-y|^2+t^2)^{{(n+s)}/{2}}}dy,$$
where  $f\ast g$ means the convolution of $f$ and $g$.   Here the constant $c(n,s)$ is the normalized constant such that $\int_{\mathbb{R}^n}p^{s}_t(x)dx=1.$
Caffarelli and  Silvestre  in \cite{Caffarelli}  proved that, with $c_{s}=\frac{\Gamma(s/2)}{2^{1-s}\Gamma(1-s/2)},$
\begin{equation}\label{eq-1.1}
(-\Delta)^{s/2}f(x)=-c_{s}\lim_{t\rightarrow{0^+}}t^{1-s}\frac{\partial u}{\partial t}(x,t).
\end{equation}
In the proof of (\ref{eq-1.1}), the key is the following identity
 \begin{equation}\label{eq2}
    \frac{1}{(2\pi)^n}\int_{\mathbb{R}^n} |\xi|^{s}|\widehat{f}(\xi)|^2d\xi
    =\frac{2^{1-s}\Gamma(1-s/2)}{\Gamma ({s}/{2})}\int_{\mathbb{R}^{n+1}_{+}}|\nabla_{(x,t)} u(x,t)|^2t^{1-s}dxdt,
    \end{equation}
which means
$$\|f\|_{\dot{H}^{s/2}}\approx \|\nabla_{(x,t)}u(\cdot,\cdot)\|_{L^{2}(\mathbb{R}^{n+1}, t^{{(1-s)}/{2}}dxdt)}.$$
To break this equivalence, in this article,  our main purpose is to identify the   space $X$ of $f$ on $\mathbb R^{n}$ such that $$\|f\|_X\lesssim \|\nabla_{(x,t)}u(\cdot,\cdot)\|_{L^{2}(\mathbb{R}^{n+1}, t^{{(1-s)}/{2}}dxdt)},$$
and to identify the space $Y$ of $u(\cdot,\cdot)$ on $\mathbb R^{n+1}_{+}$ such that
$$\|u\|_Y\lesssim\|f\|_{\dot{\Lambda}^{p,q}_\beta(\mathbb{R}^n)},$$
where $\|\cdot\|$ denotes the norms of homogeneous Besov spaces which are defined in Definition \ref{def-1.1} below.  Here and henceforth,  $A\lesssim B$ means that $A\leq CB$ for a constant $C$. $A\approx B$ means that $A\lesssim B$ and $B\lesssim A.$  Similarly, one writes ${\mathsf V}\gtrsim{\mathsf U}$ for
${\mathsf V}\ge c{\mathsf U}$ with $c$ being a constant.

Equality (\ref{eq2}) allows us to treat the fractional Sobolev/logarithmic Sobolv/Hardy inequalities as the   fractional   Sobolev/logarithmic Sobolev /Hardy trace inequalities.
For any $f\in \dot{H}^{s/2}(\mathbb{R}^n)$ and its Caffarelli-Silvestre extension
$u(x,t)=p^s_t\ast f(x) \ \forall\ (x,t)\in \mathbb{R}^{n+1}_+,$ the fractional Sobolev trace inequality  holds
\begin{equation}\label{eq3}
    \left(\int_{\mathbb{R}^n}|f(x)|^{{2n}/{(n-s)}}dx\right)^{{(n-s)}/{n}}\lesssim \int_{\mathbb{R}^{n+1}_+}|\nabla_{(x,t)} u(x,t)|^2t^{1-s}dxdt.\end{equation}
Moreover, if $\int_{\mathbb{R}^n}|f(x)|^2dx=1,$   the following fractional logarithmic Sobolev inequality holds:
\begin{equation}\label{eq4}
    \exp{ \left(\frac{s}{n}\int_{\mathbb{R}^n}|f(x)|^2\ln(|f(x)|^2)dx\right)}\lesssim \int_{\mathbb{R}^{n+1}_+}|\nabla_{(x,t)} u(x,t)|^2t^{1-s}dxdt.
    \end{equation}
There also holds the fractional Hardy inequality (or the Kato inequality) \begin{equation}\label{eq5}
    \int_{\mathbb{R}^n}\frac{|f(x)|^2}{|x|^{s}}dx\lesssim \int_{\mathbb{R}^{n+1}_{+}}|\nabla_{(x,t)} u(x,t)|^2t^{1-s}dxdt.
\end{equation}
 Xiao in \cite{Xiao2006} proved (\ref{eq3})-(\ref{eq4}) for the Poisson extension.
For the Caffarelli-Silvestre extension,  Br\"{a}ndle et al. in \cite{BC}   proved (\ref{eq3}).  Inequality (\ref{eq4}) was  proved by Nguyen in \cite{Nguyen}.  Inequalities similar to (\ref{eq5}) have been studied in  \cite{Nguyen,Beckner, Eilertsen, Herbst, Yafaev, XiaoZhai, Hajaiej Yu Zhai}.

In (\ref{eq2}), the order of the fractional derivative on left hand side   is the half  the order $s$ of the Caffarelli-Silvestre  extension. In this article, we firstly observe that \begin{equation} \label{eq2'} \int_{\mathbb{R}^{n+1}_{+}}|\nabla_{(x,t)} u(x,t)|^2t^{1-\beta}dxdt\approx \int_{\mathbb{R}^n} |\xi|^{\beta}|\widehat{f}(\xi)|^2d\xi,\quad \beta\in(0,2s),
    \end{equation}
 which allows us to  establish inequalities similar to (\ref{eq3}), (\ref{eq4})\ \&\ (\ref{eq5}) for $f\in \dot{H}^{\beta/2}(\mathbb{R}^n)$ with $\beta\in(0,2s).$
Moreover,  (\ref{eq2'}) is true when   $\nabla_{(x,t)}u(x,t)$  is replaced  by either ${\partial_{t} u(x,t)}$ or  $\nabla_{x}u(x,t)$ (more general $(-\Delta)^{-\gamma/2}$). 
  Based on these observations,  the first goal of this article is to show that
$\int_{\mathbb{R}^{n+1}_+}|\nabla_{(x,t)} u(x,t)|^2t^{1-s}dxdt$ in (\ref{eq3}), (\ref{eq4})\ \&\ (\ref{eq5}) can be   replaced by  $$\int_{\mathbb{R}^{n+1}_+}|\nabla_{(x,t)} u(x,t)|^2t^{1-\beta}dxdt,\quad
\int_{\mathbb{R}^{n+1}_+}\left| \frac{\partial u(x,t)}{\partial t}\right|^2t^{1-\beta}dxdt,\quad
\hbox{or}\quad
 \int_{\mathbb{R}^{n+1}_+}| (-\Delta)^{\gamma/2}u(x,t)|^2t^{2\gamma-\beta-1}dxdt.$$

Let $\omega:\mathbb{R}^{n+1}_+\rightarrow \mathbb{R}_+$ be a positive measurable function. Denote by $L^p(\mathbb{R}^{n+1}_+,\omega)$ the weighted Lebesgue space of all measurable functions  $g:\mathbb{R}^{n+1}_+\rightarrow \mathbb{R}$ with
$$\|g\|_{L^p(\mathbb{R}^{n+1}_+,\omega)}:=\left(\int_{\mathbb{R}^{n+1}_+}|g(x,t)|^p\omega(x,t)dxdt\right)^{1/p}<\infty.$$
In  \cite{Lombardiand}, Lombardi and Xiao  established the
affine non-sharp Poisson trace inequality for $\dot{H}^{-\beta/2}(\mathbb{R}^n)$  by using the weighted $L^p-$affine  energy
$$\mathcal{E}(g,\sigma)=c_{n,p}\left(\int_{\mathbb{S}^{n-1}}\|\nabla_{\xi}g\|^{-n}_{L^p(\mathbb{R}^{n+1}_{+},\sigma)}d\xi\right)^{-1/n},$$
where $c_{n,p}=(n\omega_{n})^{1/n}\left(\frac{n\omega_{n}\omega_{p-1}}{2\omega_{n+p-2}}\right)^{1/p}.$
We will   extend the affine non-sharp Poisson trace inequality for $\dot{H}^{-\beta/2}(\mathbb{R}^n)$  to the  Caffarelli-Silvestre extension.
The $L^p$ affine energy has been applied to establish the affine Sobolev inequalities (cf. \cite{ Lutwak-Yang-Zhang, Zhang, Zhai})
and  the affine version of the P\'{o}lya-Szeg\"{o} principle (cf. \cite{Cianchi-Lutwak-Yang-Zhang, Haberl-Schuster-Xiao}).

Let $C^{\infty}_{0}(\mathbb R^{n})$ stand for the set of all infinitely differentiable function functions with compact support.
The homogeneous Besov space  $\dot{\Lambda}^{p,q}_\beta(\mathbb{R}^n)$   are defined as  the completion of all $C^{\infty}_0(\mathbb{R}^n)$ functions with  $\|f\|_{\dot{\Lambda}^{p,q}_\beta(\mathbb{R}^n)}<\infty$, where the norm $\|f\|_{\dot{\Lambda}^{p,q}_\beta(\mathbb{R}^n)}$ is defined  as follows.
\begin{definition}\label{def-1.1}
    \item{\rm (i)}  If $(\beta,p,q)\in (0,\infty)\times (0,\infty)\times (0,\infty),$ then
    $$\|f\|_{\dot{\Lambda}^{p,q}_\beta(\mathbb{R}^n)}:=\left(\int_{\mathbb{R}^n}\|\Delta ^k_hf\|^q_{L^p(\mathbb R^{n})}|h|^{-(n+\beta q)}dh\right)^{1/q}.$$
    \item{\rm (ii)}  If $(\beta,p,q)\in (0,\infty)\times (0,\infty)\times \{\infty\},$ then
    $$\|f\|_{\dot{\Lambda}^{p,q}_\beta(\mathbb{R}^n)}:=\sup_{h\in\mathbb{R}^n\setminus\{0\}}\|\Delta^k_hf\|^q_{L^p(\mathbb R^{n})}|h|^{-\beta}.$$
 Here $k=1+[\beta], \beta=[\beta]+\{\beta\}$ with $[\beta]\in \mathbb{Z}_+,$ $\{\beta\}\in(0,1)$ and  $$\Delta^k_hf=\Delta^1_h\Delta^{k-1}_hf,\quad \Delta^1_h=f(x+h)-f(x)\quad \forall x\in\mathbb{R}^n.$$
\end{definition}

The second goal of this article is to study the Carleson type embedding for $\dot{\Lambda}^{p,q}_\beta(\mathbb{R}^n)$ via the Caffarelli-Silvestre extension. We will characterize the following embedding:
\begin{equation}\label{CarelsonEmbeddings}
    \|u\|_{L^{p_0,q_0}_{\mu}(\mathbb{R}^{n+1}_+)}\lesssim \|f\|_{\dot{\Lambda}^{p,q}_\beta(\mathbb{R}^n)}
\end{equation}
for a non-negative Borel measure/outer measure by using the Besov capacities $C^{p,q}_\beta(\cdot)$, when $\beta\in (0,1)$ or $\beta\in (0,n),$  $ p$ and $q$ satisfy one  of the following three conditions:
$$  p=q\in ({n}/{(n+\beta)},1];\quad (p,q)\in \left(1,{n}/{\beta}\right)\times (1,\infty);\quad
  (p,q)\in  \left(1,{n}/{\beta}\right)\times\{\infty\}.$$

For $0<p,q<\infty$ and a nonnegative Radon measure $\mu$ on $\mathbb{R}^{n+1}_+,$
$L^{q,p}_{\mu}(\mathbb{R}^{n+1}_+)$ and  $L^q_{\mu}(\mathbb{R}^{n+1}_+)$  denote the Lorentz space and the Lebesgue space of all functions $g(\cdot,\cdot)$ on $\mathbb{R}^{n+1}_+$, respectively, for which
$$\|g(\cdot,\cdot)\|_{L^{q,p}_{\mu}(\mathbb{R}^{n+1}_+)}:=\left\{\int_0^\infty\left(\mu\left(\{(x,t)\in \mathbb{R}^{n+1}_+:\ |g(x,t)|>\lambda \}\right)\right)^{p/q}d\lambda^p\right\}^{1/p}<\infty$$
and
$$\|g(\cdot,\cdot)\|_{L^q_{\mu}(\mathbb{R}^{n+1}_+)}:=\left(\int_{\mathbb{R}^{n+1}_+}|g(x,t)|^qd\mu\right)^{1/q}<\infty,$$
respectively. Moreover,  $L^{q,\infty}_{\mu}(\mathbb{R}^{n+1}_+)$ denotes the set of all $\mu-$measurable functions $g(\cdot,\cdot)$ on $\mathbb{R}^{n+1}_+$ with
$$\|g(\cdot,\cdot)\|_{L^{q,\infty}_{\mu}(\mathbb{R}^{n+1}_+)}:=\sup_{s>0}s\left(\mu\left(\left\{(x,t)\in \mathbb{R}^{n+1}_+:\ |g(x,t)|>s \right\}\right)\right)^{1/q}<\infty.$$

Let $C(\mathbb R^n)$ be the class of all continuous functions on $\mathbb{R}^n$. The Besov capacities $C^{p,q}_\beta(\cdot)$  have been studied in \cite{Adams1, Adams Xiao} for $q< \infty$ and in \cite{Milman Xiao} for $q=\infty.$
\begin{definition}
 \item{\rm (i)} Let $(\beta,p,q)\in (0,\infty)\times \left(0,\infty\right)\times (0,\infty).$
 For a compact set $K\subset \mathbb{R}^{n}$,  the Besov capacity $C^{p,q}_\beta(E)$ is defined as
$$C_{\beta}^{p,q}(K):=\inf\left\{\|f\|^{p}_{\dot{\Lambda}^{p,q}_\beta(\mathbb{R}^n)}:\ f\in C_0^{\infty}(\mathbb{R}^n)\text{ and } f\geq 1_K\right\}$$
and for any set $E\subset\mathbb R^{n}$, one defines
$$C_{\beta}^{p,q}(E):=\inf_{\text{ open } O\supseteq E}\sup_{\text{ compact }K\subseteq O}\left\{C^{p,q}_{\beta}(K)\right\},$$
where $1_E$ denotes the characteristic function of $E$.
 \item{\rm (ii)}  Let $(\beta,p)\in (0,1)\times \left(1,\infty\right).$  For a set $E\subset \mathbb R^{n}$, the $\infty$-Besov capacity $C^{p,\infty}_\beta(E)$ is defined as
$$C_{\beta}^{p,\infty}(E):=\inf\left\{\|f\|^{p}_{\dot{\Lambda}^{p,\infty}_\beta(\mathbb{R}^n)}:\ f\in \dot{\Lambda}^{p,\infty}_\beta(\mathbb{R}^n)\cap  C(\mathbb{R}^n) \text{ and } f\geq 1_N \text{ on a neighbourhood } N \text{ of } E\right\}.$$
\end{definition}

For any open set $O\subset\mathbb R^{n},$
the tent space $T(O)$ based on $O$ is defined as
$$T(O)=\Big\{(x,t)\in\mathbb{R}^{n+1}_+:\ B(x,t)\subset O\Big\}.$$
The variational
capacity minimizing function associated with both $C^{p,q}_\beta$ and a nonnegative
measure $\mu$ on $\mathbb{R}^{n+1}_+$
is  defined by
\begin{equation}\label{eq-1.2}
c^\beta_{p,q}(\mu; \lambda):=\inf\Big\{C^{p,q}_\beta(O):\ \hbox{open}\  O\subseteq\mathbb R^{n}, \mu(T(O))>\lambda\Big\}
\end{equation}
for $\lambda>0.$

The study of (\ref{CarelsonEmbeddings}) is mainly motivated by the work on similar embeddings via the classical Poisson kernel and the  classical/fractional heat kernels. Adams  and Xiao in  \cite{Adams Xiao} proved embeddings similar to (\ref{CarelsonEmbeddings}) via the classical Poisson kernel.  Xiao in \cite{Xiao20062, Xiao20063} studied  the  embeddings  of   $\dot{W}^{1,p}(\mathbb{R}^{n})$ and $\dot{\Lambda}^{1,1}_\beta(\mathbb{R}^n)$ into the Lebesgue space $L^{q}_\mu (\mathbb{R}^{n+1}_{+}),$ under $(p, q) \in (1,\infty) \times  (1,\infty),$   via the heat
kernel.  For fractional diffusion equations, motivated by Xiao \cite{Xiao20062},  Zhai in \cite{Zhai1}  explored the embeddings of the homogeneous Sobolev space $\dot{W}^{\beta,p}(\mathbb{R}^{n})$ into the Lebesgue space $L^{q}_\mu (\mathbb{R}^{n+1}_{+}).$
 Xiao and Zhai in \cite{XiaoZhai} explored the embeddings for the homogeneous Besov  space $\dot{\Lambda}^{p,p}_\beta(\mathbb{R}^{n})$ with $p<1$ via fractional diffusion equations.  Li et al.  in \cite{Li} studied the embedding of Sobolev spaces $\dot{W}^{\beta,p}(\mathbb{R}^{n})$ into  Lebesgue spaces $L^{q}_\mu (\mathbb{R}^{n+1}_{+})$ via the Caffarelli-Silvestre extension.

This article will be organized as follows. In Section \ref{sec-2}, we will provide some preliminaries and basic results. Section \ref{sec-3} is devoted to establish the fractional anisotropic Sobolev trace inequalities and  the affine   trace inequality for the Caffarelli-Silvestre extension.
  In Section \ref{sec-4}, we will  characterize the non-negative Random measure $\mu$ on $\mathbb{R}^{n+1}_+$ such that $ \|u(\cdot,\cdot)\|_{L^{q_0,p_0}_\mu(\mathbb{R}^{n+1}_+)}\lesssim \|f\|_{\dot{\Lambda}^{p,q}_\beta(\mathbb{R}^n)}$ in three different cases.

\section{Preliminaries and Basic Lemmas}\label{sec-2}
At first, we want to find the Fourier transform of the  fractional Poisson kernel $p_t^s(\cdot).$

\begin{proposition}\label{Prop1}
The fractional Poisson kernel $p_t^s(\cdot)$ satisfies the following properties.

    \item{\rm (i)} The Fourier transform $\widehat{p}_t^s$ can be represented as
    \begin{equation}\label{sec2-1}
    \widehat{p}_t^s(\xi)=C_{n,s} G_s(t|\xi|),
\end{equation}
with a positive constant $C_{n,s}$ and    $$G_{s}(t):=\int_0^\infty \lambda^{s/2-1}e^{-\lambda-t^2/(4\lambda)}d\lambda,\quad  t\in(0,\infty),$$
which satisfies
\begin{equation}\label{eq-2.1}
\left\{\begin{aligned}
\int_0^\infty G_s(t)^2t^{a}dt<\infty,&\quad a>-(s+1);\\
\int_0^\infty G'_s(t)^2t^{a}dt<\infty,&\quad   a>1-2s.
\end{aligned}\right.
\end{equation}

    \item{\rm (ii)}  If $a> -(s+1),$ there exists a constant $C(n,s,a)$ such that
    \begin{equation}\label{sec2-3}
    \int_0^\infty|\widehat{p}_t^s(\xi)|^2t^{a}dt
    =C(n,s,a) |\xi|^{-(a+1)}.
\end{equation}

\end{proposition}
\begin{proof}
For (i), obviously, $p^{s}_t(xt)=c(n,s)t^{-n}(|1+|x|^2)^{-{(n+s)}/{2}}$. Denote $\tau=n+s.$
For $A>0,$ it follows from the Gamma-function identity
$$A^{-\tau/2}=\frac{1}{\Gamma(\tau/2)}\int_{0}^\infty e^{-\lambda A} \lambda^{\tau/2-1}d\lambda$$
that
$$(1+|x|^2)^{-\tau/2}=\frac{1}{\Gamma(\tau/2)}\int_{0}^\infty e^{-\lambda-\lambda|x|^2} \lambda^{\tau/2-1}d\lambda.$$
Thus, we obtain
\begin{eqnarray*}
\widehat{p}_t^s(\xi)&=&t^{-n}\widehat{\delta^{t^{-1}}p^{s}_1}(\xi)\\
&=&\delta^{t}(\widehat{p^{s}_1}(\xi))\\
&=&
\frac{c(n,s)}{\Gamma(\tau/2)}\delta^t\left(\int_{0}^\infty e^{-\lambda}\widehat{(e^{-\lambda|x|^2})}(\xi) \lambda^{\tau /2-1}d\lambda\right)\\
&=&C_{n,s}\delta^t\left(\int_{0}^\infty  \lambda^{-n/2}e^{-\lambda-{|\xi|^2}/{(4\lambda)}} \lambda^{\tau/2-1}d\lambda\right)\\
&=&C_{n,s}\int_{0}^\infty e^{-\lambda-{|t\xi|^2}/{(4\lambda)}} \lambda^{s/2-1}d\lambda\\
&=&C_{n,s} G_s(t|\xi|),
\end{eqnarray*}
which gives us (\ref{sec2-1}).  See also \cite[Proposition 7.6]{Hao} for a similar result for Bessel potentials.

It follows from  \cite[page 132,  (30)]{Stein}  and  \cite[Proposition 7.6]{Hao} that there holds
\begin{equation}\label{eq-2.2}
G_{s}(t)\lesssim\left\{\begin{aligned}
&1+t^{s},&\ t\rightarrow 0;\\
&e^{-ct},&\ t\rightarrow \infty.
\end{aligned}\right.
\end{equation}

Now we verify (\ref{eq-2.1}).  Since $s+a>-1$ and $s\in (0,2),$ (\ref{eq-2.2}) implies
\begin{eqnarray*}
\int_{0}^{\infty} |G_s(t)|^{2}t^{a}dt \lesssim \int_{0}^{1}\left(1+t^{s+a}+t^{2s+a}\right)dt+\int^{\infty}_{1}e^{-2ct}t^{a}dt<\infty.
\end{eqnarray*}
Below we estimate the derivative $G'_{s}(\cdot)$. Notice that
$$G'_{s}(t)=-\frac{t}{2}\int^{\infty}_{0}\lambda^{s/2-2}e^{-\lambda-t^{2}/(4\lambda)}d\lambda.$$
Then for $t$ near $0$, letting $4\lambda=1/w$, we obtain
\begin{eqnarray*}
|t^{1-s}G'_{s}(t)|&=&\left|-\frac{t^{2-s}}{2}\int^{\infty}_{0}\lambda^{s/2-2}e^{-\lambda-t^{2}/(4\lambda)}d\lambda\right|\\
&\lesssim&t^{2-s}\int^{\infty}_{0}w^{-s/2}e^{-t^{2}w}dw,
\end{eqnarray*}
which, together with  the change of variable: $t^{2}w=v$,  gives
$$|t^{1-s}G'_{s}(t)|\lesssim \int^{\infty}_{0}v^{-s/2}e^{-v}dv<\infty,$$
where we have used the fact that $s\in(0,2)$.

Then we assume that $t\gg1$. For $N>1$, taking $2l=N+s-1$, we can apply  the change of variables $4\lambda=1/w$ and $t^{2}w=v$ again to get
\begin{eqnarray*}
|t^{N}G'_{s}(t)|&=&\left|-\frac{t^{N+1}}{2}\int^{\infty}_{0}\lambda^{s/2-2}e^{-\lambda-t^{2}/(4\lambda)}d\lambda\right|\\
&\lesssim&t^{N+1}\int^{\infty}_{0}\lambda^{s/2-2-l}e^{-t^{2}/(4\lambda)}d\lambda\\
&\lesssim&t^{N+1}\int^{\infty}_{0}w^{l-s/2}e^{-t^{2}w}dw\\
&\lesssim&\int^{\infty}_{0}v^{l-s/2}e^{-v}dv<\infty.
\end{eqnarray*}

Therefore,  we  have proved that
$$G'_{s}(t)\lesssim\left\{\begin{aligned}
t^{s-1},&\quad t\rightarrow0;\\
t^{-N},&\quad t\rightarrow \infty\ \&\ \forall N>1,
\end{aligned}\right.$$
which indicates
$$\int_{0}^{\infty} |G'_s(t)|^{2}t^{a}dt<\infty, \quad a>1-2s.$$

(ii)  For $s+a> -1,$ using  (\ref{eq-2.1}),
we have
\begin{eqnarray*}
\int_0^\infty(\widehat{p}^{s}_t(\xi))^2t^a dt&= &C_{n,s}    \int_0^\infty(G_{s}(t|\xi|))^2  (t|\xi|)^{a} d(t|\xi|)|\xi|^{-(a+1)}\\
    &= &C(n,s,a) |\xi|^{-(a+1)}.
\end{eqnarray*}
Thus (\ref{sec2-3}) holds with $C(n,s,a)=C_{n,s}\int_0^\infty G_s(t)^2t^{a}dt .$

\end{proof}

Denote  by $\mathcal{M}f$ the Hardy-Littlewood maximal operator, i.e.,
$$\mathcal{M}f(x)=\sup_{r>0}r^{-n}\int_{B(x,r)}|f(y)|dy,\quad  x\in \mathbb{R}^n.$$

\begin{lemma}\label{le-5.1}
There exists a constant $C$ depending only on $n$ and $s$ such that
$$|p^{s}_{t}\ast f(x-y)|\leq C (1+|y|^{2}/t)\mathcal{M}f(x)\quad \forall x\  \&\  y\in\mathbb R^{n}\quad \forall f\in C^\infty_0(\mathbb{R}^n).$$
\end{lemma}
\begin{proof}
For any ball $Q\subset \mathbb R^{n}$, denote by $x_{Q}$ the center of $Q$ and by $r_{Q}$ the radius of $Q$. We can see  that if $x/t\in Q$ then $|x-tx_{Q}|<tr_{Q}$. This implies that
$$\mathcal{M}f(t\cdot)(x/t)\leq\sup_{x\in tQ}\frac{1}{|tQ|}\int_{tQ}|f(u)|du\leq \mathcal{M}f(x).$$
We only need to prove
\begin{equation}\label{eq-5.5}
|p^{s}_{1}\ast f(x-y)|\lesssim (1+|y|)^{n}\mathcal{M}f(x).
\end{equation}
In fact, if (\ref{eq-5.5}) holds, then by the change of variable $z/t=u$, we obtain
\begin{eqnarray*}
|p^{s}_{t}\ast f(x-y)|&\approx&\Big|\int_{\mathbb R^{n}}\frac{t^{s}}{(t^{2}+|x-y-z|^{2})^{(n+s)/2}}f(z)dz\Big|\\
&\approx&\left|\int_{\mathbb R^{n}}\frac{1}{t^{n}}\frac{1}{(1+|\frac{x-y}{t}-u|^{2})^{(n+s)/2}}f(tu)t^{n}du\right|\\
&\lesssim&(1+|y|/t)^{n}\mathcal{M}f(t\cdot)(x/t)\\
&\lesssim&(1+|y|/t)^{n}\mathcal{M}f(x).
\end{eqnarray*}
Below we prove (\ref{eq-5.5}). We first estimate the kernel $p^{s}_{1}(x-y)$. If $|x|<|y|$, it is obvious that $p^{s}_{1}(x-y)\leq 1$. On the other hand, for $|x|\geq |y|$, noting that $||x|-|y||\leq |x-y|$, it holds
$$p^{s}_{1}(x-y)\lesssim \frac{1}{(1+(|x|-|y|)^{2})^{(n+s)/2}}.$$
Set the decreasing radial majorant function of $p^{s}_{1}(x-y)$ as
$$\varphi_{y}^s(x):=\left\{\begin{aligned}
&&\frac{1}{(1+(|x|-|y|)^{2})^{(n+s)/2}},\quad |x|\geq |y|;\\
&&1,\quad  |x|<|y|.
\end{aligned}\right.$$
With a slight abuse of notation, let us write $\varphi^{s}_{y}(x)=\varphi^{s}_{y}(r)$ if $|x|=r$. We can get
\begin{eqnarray*}
|p^{s}_{1}\ast f(x-y)|&\lesssim&\int_{\mathbb R^{n}}\frac{1}{(1+|x-y-z|^{2})^{(n+s)/2}}|f(z)|dz\\
&\lesssim&\int_{\mathbb R^{n}}\varphi_{y}^s(x-z)|f(z)|dz\\
&\approx&\sum^{\infty}_{k=-\infty}\int_{2^{k}<|x-z|\leq 2^{k+1}}\varphi_{y}^s(x-z)|f(z)|dz\\
&\lesssim&\sum^{\infty}_{k=-\infty}\varphi_{y}^s(2^{k})\int_{|x-z|<2^{k+1}}|f(z)|dz\\
&\lesssim&\mathcal{M} f(x)\sum^{\infty}_{k=-\infty}\varphi_{y}^s(2^{k})2^{n(k+1)}\\
&\lesssim&\mathcal{M}f(x)\sum^{\infty}_{k=-\infty}\varphi_{y}^s(2^{k})\int^{2^{k}}_{2^{k-1}}r^{n-1}dr\\
&\lesssim&\mathcal{M}f(x)\|\varphi_{y}^s\|_{L^1(\mathbb R^{n})}.
\end{eqnarray*}
By a direct computation, we obtain
\begin{eqnarray*}
\|\varphi_{y}^s\|_{1}&\approx&\int_{|x|<|y|}1dx+\int_{|x|\geq|y|}\frac{1}{(1+(|x|-|y|)^{2})^{(n+s)/2}}dx\\
&\lesssim&|y|^{n}+\int_{|x|\geq |y|}\frac{1}{(1+|x|^{2}+|y|^{2}-2|x||y|)^{(n+s)/2}}dx\\
&\lesssim&|y|^{n}+\int^{\infty}_{0}\frac{(r+|y|)^{n-1}dr}{(1+r^{2})^{(n+s)/2}}\\
&\lesssim&|y|^{n}+\sum^{n-1}_{k=0}C^{k}_{n-1}|y|^{n-1-k}\int^{\infty}_{0}\frac{r^{k}}{(1+r^{2})^{(n+s)/2}}dr\\
&\lesssim&|y|^{n}+\sum^{n-1}_{k=0}C^{k}_{n-1}|y|^{n-1-k}\\
&\leq& C(1+|y|)^{n}.
\end{eqnarray*}

\end{proof}

 In \cite{Xiao2006}, for the heat kernel and variational capacities $cap_p(\cdot)$,  Xiao established some basic properties concerning tents, non-tangential maximal functions and $cap_p(\cdot)$. Following the idea of \cite[Lemma 2.4]{Xiao2006}, we can deduce  the following properties for the fractional Poisson kernel $p^{s}_{t}(\cdot)$ and the related Besov capacities $C^{p,q}_\beta(\cdot)$. Parts of the following result were proved for Sobolev functions in \cite{Li}.

\begin{lemma}\label{lemma 2}
Let $s \in(0,2) $ and $\beta\in (0,n).$ Given $f\in C^\infty_0(\mathbb{R}^n)$, denote the Caffarelli-Silvestre extension of $f$ by
$u(x,t)=p_t^{s}\ast f(x).$  For $\lambda>0$ and a non-negative measure $\mu$ on $\mathbb R^{n+1}_{+},$ let
$$E_{\lambda,s}(f)=\left\{(x,t)\in \mathbb{R}^{n+1}_+:\ |u(x,t)|>\lambda\right\}$$
and
$$O_{\lambda,s}(f)=\Big\{x\in \mathbb{R}^{n}:\ \sup_{|y-x|<t}|u(y,t)|>\lambda\Big\}.$$
Then the following four statements are true.
\item{\rm(i)}   $\mu(E_{\lambda,s}(f))\leq \mu(T(O_{\lambda,s}(f))).$
\item{\rm (ii)}
For any natural number $k$,
    $$\mu\left(E_{\lambda,s}(f) \cap T(B(0,k))\right)\leq \mu\left(T\left(O_{\lambda,s}(f)\cap B(0,k)\right)\right).$$
\item{\rm (iii)} Let $c^\beta_{p,q}$ be the variational
capacity minimizing function defined by (\ref{eq-1.2}). For any $p,q\in (0,\infty)$  and any natural number $k,$
    $$C^{p,q}_\beta\left(O_{\lambda,s}(f)\cap B(0,k)\right)\geq c^\beta_{p,q}(\mu; \mu\left(T\left(O_{\lambda,s}(f)\cap B(0,k)\right)\right).$$
\item{\rm (iv)} There exists a  constant $\theta_{n,s}>0$ such that
    $$\sup_{|y-x|<t}|u(y,t)|\leq \theta_{n,s} \mathcal{M}f(x),\quad  x\in \mathbb{R}^n.$$

\item{\rm (v)} There exists a  constant $\eta_{n,s}>0$ such that if $(x,t)\in T(O)$, then
$$(p^{s}_t\ast |f|)(x) \geq \eta_{n,s},$$
where $O$ is a bounded open set contained in $\text{ Int }(\{x\in \mathbb{R}^n: f(x)\geq 1\}).$

\end{lemma}
\begin{proof}
    (i) and (ii).
Let $\Gamma(x)=\Big\{(y,t)\in\mathbb{R}^{n+1}_{+}:\ |y-x|<t\Big\}$. The non-tangential maximal function related with $p^{s}_{t}$ is defined as
$$N(f)(x):=\sup_{|y-x|<t}|p^{s}_{t}\ast f(y)|.$$
Then $O_{\lambda,s}=\Big\{x\in\mathbb R^{n}:\ N(f)(x)>\lambda\Big\}.$ For any $x_{0}\in O_{\lambda,s}$, $N(f)(x_{0})>\lambda$, i.e., there exists a $y_{0}$ such that $|y_{0}-x_{0}|<t$ and $|p^{s}_{t}\ast f|(y_{0})>\lambda$. For arbitrary $x\in B(y_{0}, t)$, it holds
$$N(f)(x):=\sup_{(z,t)\in\Gamma(x)}|p^{s}_{t}\ast f(z)|\geq |p^{s}_{t}\ast f(y_{0})|>\lambda,$$
which indicates that $B(y_{0}, t)\in O_{\lambda,s}$ and $N(f)$ is lower semi-continuous. We can see that $O_{\lambda,s}$ is an open set in $\mathbb R^{n}$.
    So, we get
    $$\left\{\begin{aligned}
    &E_{\lambda,s}(f)\subseteq T(O_{\lambda,s}(f));\\
     &\mu(E_{\lambda,s}(f))\leq \mu(T(O_{\lambda,s}(f))),
    \end{aligned}\right.$$
    and
    \begin{eqnarray*}
    \mu(E_{\lambda,s}(f)\cap T(B(0,k)))\leq \mu(T(O_{\lambda,s}(f))\cap T(B(0,k)))=\mu(T(O_{\lambda,s}(f)\cap B(0,k))).
    \end{eqnarray*}

   (iii) It follows from the definition of $c^\beta_{p,q}(\mu;t).$

    (iv)  By Lemma \ref{le-5.1}, for $y$ satisfying $|x-y|<t$,
$$p^{s}_{t}\ast f(x-y)\lesssim (1+|y|/t)^{n}\mathcal{M}f(x),$$
which gives
$$N(f)(x)=\sup_{|y|<t}|p^{s}_{t}\ast f(x-y)|\lesssim \mathcal{M}f(x)\sup_{|y|<t}(1+|y|/t)^{n}\lesssim 2^{n}\mathcal{M}f(x).$$

    (v)  For any $(x,t)\in T(O)$, we have
$$B(x,t)\subseteq O\subset \text{Int}(\{x\in\mathbb R^{n}: f(x)\geq1\}).$$
We can see that for $|x|\leq \sigma t$,
$$p^{s}_{t}(x)\approx\frac{t^{s}}{(t^{2}+|x|^{2})^{(n+s)/2}}\gtrsim\frac{t^{s}}{(t^{2}+(\sigma t)^{2})^{(n+s)/2}}\gtrsim\frac{1}{t^{n}}.$$
Then
\begin{eqnarray*}
p^{s}_{t}\ast |f|(x)&=&\int_{\mathbb R^{n}}p^{s}_{t}(x-y)|f|(y)dy\\
&\gtrsim&\frac{1}{t^{n}}\int_{B(x,\sigma t)\cap \text{ Int }(\{x\in\mathbb R^{n}:\ f(x)\geq1\})}|f(y)|dy\\
&\gtrsim&\frac{1}{t^{n}}\Big|B(x,\sigma t)\cap \text{ Int }\Big(\Big\{x\in\mathbb R^{n}:\ f(x)\geq1\Big\}\Big)\Big|.
\end{eqnarray*}

If $\sigma>1$, then
$$B(x,t)=B(x, t)\cap\text{ Int }(\{x\in\mathbb R^{n}:\ f(x)\geq1\})\subset B(x,\sigma t)\cap\text{ Int }(\{x\in\mathbb R^{n}:\ f(x)\geq1\}).$$

If $\sigma\leq1$, then
$$B(x,\sigma t)\cap\text{ Int }(\{x\in\mathbb R^{n}:\ f(x)\geq1\})=B(x,\sigma t).$$
Hence $p^{s}_{t}\ast |f|(x)\geq \eta_{n,s}$ for a constant $\eta_{n,s}$.

    \end{proof}

    We need the following week/strong type estimates for Besov capacities $C^{p,q}_{\beta}(\cdot)$.    Denote
    $$p\lor q=\max\{p,q\}\text{ and } p\land q=\min\{p,q\}.$$
    \begin{lemma}\label{lemma 24}
    The following results hold.
    \item{\rm (i)} When $\beta\in (0,1)\ \&  \ p=q\in\left({n}/{(n+\beta)},1\right),$ or $(\beta,p,q)\in(0,n)\times \left(1,{n}/{\beta}\right)\times (1,\infty),$ $$\int_{0}^\infty\left(C^{p,q}_\beta\left(\left\{x\in\mathbb{R}^n:\ |\mathcal{M}f(x)|> \lambda\right\}\right)\right)^{1\lor ({q}/{p})}d\lambda^{p\lor q}\lesssim \|f\|^{p\lor q}_{\dot{\Lambda}^{p,q}_\beta(\mathbb{R}^n)}
    \quad\forall\ f\in C^{\infty}_0(\mathbb{R}^n).$$

    \item{\rm (ii)} When $(\beta, p,q)\in (0,1)\times \left(1,{n}/{\beta}\right)\times \{\infty\},$ $$\sup_{\lambda\in(0,\infty)}\lambda^p
C^{p,q}_\beta\left(\left\{(x,t)\in \mathbb{R}^{n}:\ |\mathcal{M}f(x)|>\lambda\right\}\right)\lesssim \|f\|^p_{\dot{\Lambda}^{p,q}_{\beta}(\mathbb{R}^n)}\quad \forall\ f\in \dot{\Lambda}^{p,q}_{\beta}(\mathbb{R}^n)\cap C(\mathbb{R}^n).$$

    \end{lemma}

\begin{proof}
Statement (i) is due to Maz'ya \cite{Mazya} when $p=q>1.$ When $1\leq p\leq q<\infty$ and $0<\beta<1,$ Wu \cite{Wu} proved (i).  Adams and Xiao \cite{Adams Xiao} established (i) when $0<\beta<\infty$ and $(p,q)\in (1, n/\beta)\times (1,\infty).$ Xiao and Zhai
\cite{Xiao Zhai 2} showed that (i) holds when $0<\beta<1$ and ${n}/{(n+\beta)}<p=q<1.$

For (ii), it follows from \cite[Proposition 2.8]{Milman Xiao}  that for $\lambda>0$ and  $f\in \dot{\Lambda}^{p,q}_{\beta}(\mathbb{R}^n)\cap C(\mathbb{R}^n),$
$$\left(\lambda^pC^{p,q}_\beta\left(\left\{(x,t)\in \mathbb{R}^{n}:\ |\mathcal{M}f(x)|>\lambda\right\}\right)\right)^{1/p}\lesssim \|\mathcal{M}f\|_{\dot{\Lambda}^{p,q}_{\beta}(\mathbb{R}^n)}
\lesssim\|f\|_{\dot{\Lambda}^{p,q}_{\beta}(\mathbb{R}^n)},$$
where in the last inequality we have used the fact that the maximal function is bounded on $L^p(\mathbb R^{n})$.
\end{proof}

\section{Fractional Trace Inequalities via the Caffarelli-Silvestre Extension}\label{sec-3}
\subsection{Fractional trace inequalities involving $\nabla_{(x,t)} u(x,t)$}\label{sec-3.1}
\begin{theorem}\label{them4}
Let $f\in \dot{H}^{\beta/2}(\mathbb{R}^n)$ with $\beta\in (0,\min\{n,2s\}).$ Denote by
$u(x,t)=p_t^{s}\ast f(x)$ the Caffarelli-Silvestre extension of $f$.
    \item{\rm (i)}
There holds
\begin{equation}\label{eq5'}
\left(\int_{\mathbb{R}^n}|f(x)|^{{2n}/{(n-\beta)}}dx\right)^{1-\beta/n}\lesssim \int_{\mathbb{R}^{n+1}_+}|\nabla_{(x,t)} u(x,t)|^2t^{1-\beta}dxdt.
\end{equation}
\item{\rm (ii)}  If $\|f\|_{L^2(\mathbb{R}^n)}=1,$ there holds
\begin{equation}\label{eq6'}
\exp{ \left(\frac{\beta}{n}\int_{\mathbb{R}^n}|f(x)|^2\ln(|f(x)|^2)dx\right)} \lesssim\int_{\mathbb{R}^{n+1}_+}|\nabla_{(x,t)} u(x,t)|^2t^{1-\beta}dxdt.
\end{equation}
\item{\rm (iii)} There  holds
\begin{equation}\label{eq7'}
\int_{\mathbb{R}^n}|f(x)|^2\frac{dx}{|x|^{\beta}}\lesssim \int_{\mathbb{R}^{n+1}_+}|\nabla_{(x,t)} u(x,t)|^2t^{1-\beta}dxdt.
\end{equation}
\end{theorem}

\begin{proof}
In order to prove Theorem \ref{them4}, we need to establish the following result.
For $\beta\in(0,2s),$ there holds
\begin{equation}\label{eq12'}
    \int_{\mathbb{R}^{n+1}_+}|\nabla_{(x,t)} u(x,t)|^2t^{1-\beta}dxdt\approx \int_{\mathbb{R}^n}|\xi|^{\beta}(\widehat{f}(\xi))^2d\xi.
\end{equation}
In fact,
since $u(x,t)=p^{s}_t\ast f(x)$, (\ref{eq-2.1}) of  Proposition \ref{Prop1} implies
\begin{eqnarray*}
    &&\int_{\mathbb{R}^{n+1}_+}  |\nabla_{(x,t)} u(x,t)|^2t^{1-\beta}dxdt\\
    &&=
    \int_{\mathbb{R}^{n+1}_+}  |\nabla_x u(x,t)|^2t^{1-\beta}dxdt+\int_{\mathbb{R}^{n+1}_+}  \left|\frac{ \partial u(x,t)}{\partial t}\right|^2t^{1-\beta}dxdt\\
&&=  \int_0^\infty\int_{\mathbb{R}^{n}}|\xi|^2|\widehat{u}(\xi,t)|^2d\xi t^{1-\beta} dt+ \int_0^\infty\int_{\mathbb{R}^{n}}\left|\frac{\partial \widehat{u}(\xi,t)}{\partial t}\right|^2d\xi t^{1-\beta} dt
    \\
    &&=   \int_0^\infty\int_{\mathbb{R}^{n}}|\xi|^2(\widehat{p}^{s}_t(\xi))^2 (\widehat{f}(\xi))^2d\xi t^{1-\beta} dt+    \int_0^\infty\int_{\mathbb{R}^{n}}\left(\frac{\partial(\widehat{p}^{s}_t(\xi)}{\partial t}\right)^2 (\widehat{f}(\xi))^2d\xi t^{1-\beta} dt\\
    &&\approx    \int_0^\infty\int_{\mathbb{R}^n}|\xi|^2(G_s(t|\xi|))^2 (\widehat{f}(\xi))^2d\xi t^{1-\beta} dt+    \int_0^\infty\int_{\mathbb{R}^{n}}\left(\frac{\partial(G_{s}(t|\xi|)}{\partial t}\right)^2 (\widehat{f}(\xi))^2d\xi t^{1-\beta} dt\\
            &&\approx
    \int_{\mathbb{R}^{n}}\left(\int_0^\infty(G_{s}(u))^2  u^{1-\beta} du+\int_0^\infty(G'_{s}(u))^2  u^{1-\beta} du\right)|\xi|^{\beta}(\widehat{f}(\xi))^2d\xi\\
    &&\approx \int_{\mathbb{R}^{n}}|\xi|^{\beta}(\widehat{f}(\xi))^2d\xi.
\end{eqnarray*}
In the last equality, we have used the requirement $s+1-\beta>-1$ and $1-\beta>1-2s$ which means that $\beta<2+s$ and $\beta<2s.$ Since $s\in(0,2),$ $2s<2+s$ and $\beta\in (0,2s).$
It follows from (\ref{eq12'}) and the well-known fractional Sobolev inequality
$$\|f\|_{L^{{2n}/{(n-\beta)}}(\mathbb{R}^n)}\leq B(n,\beta)\|(-\Delta)^{\beta/4}f\|_{L^2(\mathbb{R}^n)}$$
for $\beta\in(0,n)$ and some constant $B(n,\beta)$ that (\ref{eq5'}) holds.

Now, we want to prove (\ref{eq6'}). Let $p={n(r-2)}/{\beta}, 2<r<{2n}/{(n-\beta)}$ and $\beta\in (0,\min\{n,2s\}).$ Then the H\"{o}lder inequality implies
$$\|f\|_{L^r(\mathbb R^{n})}^r=\int_{\mathbb{R}^n}|f(x)|^{p}|f(x)|^{r-p}dx\leq \|f\|^p_{L^{{2n}/{(n-\beta)}}(\mathbb R^{n})}\left(\int_{\mathbb{R}^n}|f(x)|^2dx\right)^{1-{p(n-\beta)}/{(2n)}}.$$
If $\|f\|_{L^2(\mathbb{R}^n)}=1$, then
$$\left(\int_{\mathbb{R}^n}|f(x)|^{r-2}|f(x)|^2dx\right)^{{1}/{(r-2)}}=\left(\int_{\mathbb{R}^n} |f(x)|^rdx\right)^{{1}/{(r-2)}}\leq\|f\|_{L^{{2n}/{(n-\beta)}}(\mathbb R^{n})}^{{n}/{\beta}}.$$
So, the inequality (\ref{eq5'}) implies that for a positive  constant $A(n,s,\beta),$
$$\left(\int_{\mathbb{R}^n}|f(x)|^{r-2}|f(x)|^2dx\right)^{{1}/{(r-2)}}\leq \left(A(n,s,\beta) \int_{\mathbb{R}^{n+1}_+}|\nabla_{(x,t)} u(x,t)|^2t^{1-\beta}dxdt.
\right)^{{n}/{(2\beta)}}.$$
 Since $\|f\|_{L^2(\mathbb{R}^n)}=1,$ $(f(x))^2dx$ can be treated as a probability measure on $\mathbb{R}^n.$ Thus, (\ref{eq6'}) can be obtained by letting $r\rightarrow 2 $ in the previous inequality.

Inequality (\ref{eq7'}) follows from (\ref{eq12'})  and the fractional Hardy inequality
$$\left\|\frac{f(\cdot)}{|\cdot|^{\beta/2}}\right\|_{L^2(\mathbb{R}^n)}\leq H \|(-\Delta)^{\beta/4}f\|_{L^2(\mathbb{R}^n)}$$
which is a special case of \cite[(3.1) in Theorem 3.1]{XiaoZhai}.

\end{proof}

\subsection{Fractional anisotropic trace inequalities involving $\partial_{t} u(x,t)$}\label{sec-3.2}

\begin{theorem}\label{them1}
Let $f\in \dot{H}^{\beta/2}(\mathbb{R}^n)$ with $\beta\in (0,\min\{n,2s\}).$ Denote by
$u(x,t)=p_t^{s}\ast f(x)$ the Caffarelli-Silvestre extension of $f$.
    \item{\rm (i)}
There holds
\begin{equation}\label{eq8'}
\left(\int_{\mathbb{R}^n}|f(x)|^{{2n}/{(n-\beta)}}dx\right)^{1-{\beta}/{n}}\lesssim  \int_{\mathbb{R}^{n+1}_+}\left| \frac{\partial u(x,t)}{\partial t}\right|^2t^{1-\beta}dxdt.
\end{equation}
\item{\rm (ii)}  If $\|f\|_{L^2(\mathbb{R}^n)}=1,$ there holds
\begin{equation}\label{eq9'}
\exp{ \left(\frac{\beta}{n}\int_{\mathbb{R}^n}|f(x)|^2\ln(|f(x)|^2)dx\right)} \lesssim \int_{\mathbb{R}^{n+1}_+}\left| \frac{\partial u(x,t)}{\partial t}\right|^2t^{1-\beta}dxdt.
\end{equation}
\item{\rm (iii)} There holds
\begin{equation}\label{eq 10'}
\int_{\mathbb{R}^n}\frac{|f(x)|^2}{|x|^{\beta}}dx\lesssim \int_{\mathbb{R}^{n+1}_+}\left| \frac{\partial u(x,t)}{\partial t}\right|^2t^{1-\beta}dxdt.
\end{equation}
\end{theorem}
\begin{proof}
In order to prove Theorem \ref{them1}, we need to establish the following result.
For $\beta\in(0,\min\{n,2s\}),$ there exists $a(n,s,\beta)$ such that
\begin{equation}\label{eq12}
    \int_{\mathbb{R}^{n+1}_+}\left|\frac{\partial u(x,t)}{\partial t}\right|^2t^{1-\beta}dxdt=a(n,s,\beta)\int_{\mathbb{R}^n}| \xi|^{\beta}(\widehat{f}(\xi))^2d\xi.
\end{equation}

In fact,
noting that $u(x,t)=p^{s}_t\ast f(x)$, we can apply Proposition \ref{Prop1} to deduce that
\begin{eqnarray*}
    \int_{\mathbb{R}^{n+1}_+}\left|\frac{\partial u(x,t)}{\partial t}\right|^2t^{1-\beta}dxdt&=&
\int_{\mathbb{R}^{n+1}_+}  \left|\frac{ \partial u(x,t)}{\partial t}\right|^2t^{1-\beta}dxdt\\
&=&    \int_0^\infty\int_{\mathbb{R}^{n}}\left|\frac{\partial \widehat{u}(\xi,t)}{\partial t}\right|^2d\xi t^{1-\beta} dt
    \\
    &=&    \int_0^\infty\int_{\mathbb{R}^{n}}\left(\frac{\partial(\widehat{p}^{s}_t(\xi)}{\partial t}\right)^2 (\widehat{f}(\xi))^2d\xi t^{1-\beta} dt\\
            &=&
    \int_{\mathbb{R}^{n}}\left(\int_0^\infty(G'_{s}(u))^2  u^{1-\beta} du\right)|\xi|^{\beta}(\widehat{f}(\xi))^2d\xi\\
    &= & a(n,s,\beta) \int_{\mathbb{R}^{n}}|\xi|^{\beta}(\widehat{f}(\xi))^2d\xi.
\end{eqnarray*}
Then, Theorem \ref{them1} can be proved in a way similar to that of Theorem \ref{them4}.
\end{proof}

\subsection{Fractional anisotropic  trace inequalities involving $(-\Delta)^{\gamma/2}u(x,t)$}\label{sec-3.3}

\begin{theorem}\label{them2}
For   $f\in\dot{H}^{\beta/2}(\mathbb{R}^n)$ with $\beta\in (0,n),$ $\gamma>\max\left\{{(\beta-s)}/{2},0\right\}.$ Denote by
$u(x,t)=p_t^{s}\ast f(x)$ the Caffarelli-Silvestre extension of $f$.
 \item{\rm (i)}
    There holds
\begin{equation}\label{7}
\left(\int_{\mathbb{R}^n}|f(x)|^{{2n}/{(n-\beta)}}dx\right)^{1-\beta/n}\lesssim \int_{\mathbb{R}^{n+1}_+}| (-\Delta)^{\gamma/2}u(x,t)|^2t^{2\gamma-\beta-1}dxdt.
\end{equation}
\item{\rm (ii)}   If $\|f\|_{L^2(\mathbb{R}^n)}=1,$ there holds
\begin{equation}\label{8}
\exp{ \left(\frac{\beta}{n}\int_{\mathbb{R}^n}|f(x)|^2\ln(|f(x)|^2)dx\right)}\lesssim\int_{\mathbb{R}^{n+1}_+}| (-\Delta)^{\gamma/2}u(x,t)|^2t^{2\gamma-\beta-1}dxdt .
\end{equation}
\item{\rm (iii)} There holds
\begin{equation}\label{10'}
\int_{\mathbb{R}^n}\frac{|f(x)|^2}{|x|^{\beta}}dx\lesssim \int_{\mathbb{R}^{n+1}_+}| (-\Delta)^{\gamma/2}u(x,t)|^2t^{2\gamma-\beta-1}dxdt.
\end{equation}
\end{theorem}
\begin{proof}
 We first establish the following result: for $\beta\in(0,n),$ and $\gamma>\max\left\{{(\beta-s)}/{2},0\right\},$  there holds \begin{equation}\label{eq15}
    \int_{\mathbb{R}^{n+1}_+}|(-\Delta)^{\gamma/2}u(x,t)|^2t^{2\gamma-\beta-1}dxdt\approx \int_{\mathbb{R}^n}|\xi|^{\beta}(\widehat{f}(\xi))^2d\xi.
\end{equation}

In fact, similar to the proof of  (\ref{eq12}), (i) of Proposition \ref{Prop1} implies
\begin{eqnarray*}
    \int_{\mathbb{R}^{n+1}_+}  |(-\Delta)^{\gamma/2}u(x,t)|^2t^{2\gamma-\beta-1}dxdt
    &\approx&
    \int_0^\infty\int_{\mathbb{R}^{n}}|\xi|^{2\gamma}|\widehat{u}(\xi,t)|^2d\xi t^{2\gamma-\beta-1} dt
    \\
    &\approx&
    \int_0^\infty\int_{\mathbb{R}^{n}}|\xi|^{2\gamma}(\widehat{p}^{s}_t(\xi))^2 (\widehat{f}(\xi))^2d\xi t^{2\gamma-\beta-1} dt\\
        &\approx&
    \int_0^\infty\int_{\mathbb{R}^{n}}(G_{s}(t|\xi|))^2 |\xi|^{2\gamma}(\widehat{f}(\xi))^2d\xi t^{2\gamma-\beta-1} dt\\
        &\approx&
    \int_{\mathbb{R}^{n}}\left(\int_0^\infty(G_{s}(t|\xi|))^2  (t|\xi|)^{2\gamma-\beta-1} d(t|\xi|)\right)|\xi|^{\beta}(\widehat{f}(\xi))^2d\xi\\
    &\approx &  \int_{\mathbb{R}^{n}}|\xi|^{\beta}(\widehat{f}(\xi))^2d\xi.
\end{eqnarray*}
Similar to the proof  of Theorem \ref{them4},
we can derive  Theorem \ref{them2} by using (\ref{eq15}).
\end{proof}

\subsection{Affine fractional Poisson trace inequality}
In this section, we will prove an affine fractional Poisson  trace inequality by applying the sharp affine weighted $L^p$ Sobolev inequality established in \cite{Haddad}. Let $\alpha>0$. Define a function $\sigma$ on $\mathbb R^{n+1}_{+}$ as $\sigma(x,t):=t^{\alpha} \quad \forall\ (x,t)\in \mathbb{R}^{n+1}_+$. For $1\leq p<n+\alpha+1$, denote by $L^{p}(\mathbb R^{n+1}_{+}, \sigma)$ the weighted Lebesgue space corresponding to the weight function $\sigma$.
\begin{theorem}{\rm (\cite[Theorem 1.1]{Haddad})}\label{affineLp}
Let $\alpha\geq 0$, $1\leq p<n+\alpha+1$ and $p^\star_\alpha={p(n+\alpha+1)}/{(n+\alpha+1-p)}$.
There exists a sharp constant $J(n,p,\alpha)$ such that
$$\|g(\cdot,\cdot)\|_{L^{p^\star_\alpha}(\mathbb{R}^{n+1}_+,\sigma)}\leq J(n,p,\alpha)(\mathcal{E}_p(g,\sigma))^{{n}/{(n+\alpha+1)}}\Big\|\frac{\partial g}{\partial t }(\cdot,\cdot)\Big\|^{(\alpha+1)/(n+\alpha+1)}_{L^p(\mathbb{R}^{n+1}_+,\sigma)}.$$
The equality in the above inequality holds if
\begin{equation*}
    g(x,t)=\left\{
\begin{aligned}
&\frac{c}{(1+|\gamma t|^{1+1/p}+|B(x-x_0)|^{1+1/p})^{{(1+n+\alpha-p)}/{p}} },&\  \ p>1;\\
&c1_{\mathbb{B}^{n+1}}(\gamma t,  B(x-x_0)),&\   p=1
\end{aligned}
\right.
\end{equation*}
for some quadruple $$(c,|\gamma|, x_0,B)\in \mathbb{R}\times \mathbb{R}_+\times \mathbb{R}^{n}\times GL_{n},$$
where $1_{\mathbb{B}^{n+1}}$ is the characteristic function of the unit ball in $\mathbb{R}^{n+1}$ and $GL_n$ denotes the set of all invertible real $n\times n-$matrices.

\end{theorem}

\begin{theorem}\label{them3} Let
 $p={2(n+\beta)}/{(n+\beta+2)}$ and  $\beta\geq 1.$
For any $f\in C^\infty_0(\mathbb{R}^n)$   and its Caffarelli-Silvestre extension
$u(x,t)=p^s_t\ast f(x),$
there holds,
\begin{equation}\label{equ-3.13}
\|f\|_{\dot{H}^{-\beta/2}(\mathbb{R}^n)}\lesssim
 (\mathcal{E}_{p}(u,t^{\beta-1}))
^{{n}/{(n+\beta)}}
\left\|\frac{\partial u}{\partial t}(\cdot,\cdot)\right\|^{{\beta}/{(n+\beta)}}_{L^p(\mathbb{R}^{n+1}_+,\ t^{\beta-1})}. \end{equation}
\end{theorem}
\begin{proof}
Equality (\ref{sec2-3}) implies
\begin{eqnarray*}
\|u(\cdot,\cdot)\|^2_{L^2(\mathbb{R}^{n+1}_+, t^a)}
\approx\int_{\mathbb{R}^n}|\xi|^{-(a+1)}|\widehat{f}(\xi)|^2d\xi
\end{eqnarray*}
for $a>-(s+1).$ In the following, taking $a=\beta-1\geq0$ in Theorem \ref{affineLp}, we can get $\sigma:=t^{a}=t^{\beta-1}$ and
\begin{eqnarray*}
\|f\|_{\dot{H}^{-\beta/2}(\mathbb{R}^n)}&=&\|(-\Delta)^{{-\beta}/{4}}f\|_{L^2(\mathbb{R}^n)}\\
&\approx&\|u(\cdot,\cdot)\|_{L^2(\mathbb{R}^{n+1}_+, \sigma)}\\
&\lesssim& (\mathcal{E}_p(u,\sigma))^{{n}/{(n+\beta)}}\left\|\frac{\partial u}{\partial t}(\cdot,\cdot)\right\|^{{\beta}/{(n+\beta)}}_{L^p(\mathbb{R}^{n+1}_+,\sigma)}.
\end{eqnarray*}
\end{proof}

\subsection{Remarks on the general case $p \geq 1 $ }
In Theorems \ref{them4}\ \&\ \ref{them1}\ \&\ \ref{them2},  the scope of $(p,\beta)$ is $p=2$ and $0<\beta<\min\{n,2s\}$. We can generalize the fractional trace type  inequalities  Theorems \ref{them4}\ \&\ \ref{them1}\ \&\ \ref{them2} to the general index $p\in[1,\infty), $ $ \beta\in (0, 2n/p).$ Moreover, a general affine fractional Poisson trace inequality similar to   Theorem \ref{them3} holds  for $f\in \dot{\Lambda}^{p^{\ast},p^{\ast}}_{-\beta/2}(\mathbb{R}^n)$, where
$p^{\ast}$ satisfies $p^{\ast}\geq \max\{2/\beta,1\}$ and $1/p=1/p^{\ast}+2/(2n+p^{\ast}\beta)$.

Denote by
$u(x,t)=p_t^{s}\ast f(x)$ the
Caffarelli-Silvestre extension of $f$. Applying \cite[Theorems 1.1 \& 1.3]{BuiCandy},  Lenzmann-Schikorra in \cite[Proposition 10.8]{LS} obtained the following characterizations of Besov spaces $\dot{\Lambda}^{p,p}_{\beta}(\mathbb R^{n})$ with $p\in (1,\infty),$
\begin{equation}\label{eq-3.1}
\left(\int_{\mathbb R^{n+1}_{+}}|\nabla_{x}u(x,t)|^{p}t^{p-1-p\beta/2}dxdt\right)^{1/p}\approx \|f\|_{\dot{\Lambda}^{p,p}_{\beta/2}(\mathbb R^{n})}, \quad\hbox{when}\quad  \beta\in (0,2);
\end{equation}
\begin{equation}\label{eq-3.2}
\left(\int_{\mathbb R^{n+1}_{+}}|\partial_{t}u(x,t)|^{p}t^{p-1-p\beta/2}dxdt\right)^{1/p}\approx \|f\|_{\dot{\Lambda}^{p,p}_{\beta/2}(\mathbb R^{n})}, \quad\hbox{when}\quad  \beta\in(0,2s);
\end{equation}
and
\begin{equation}\label{equ-3.15}
\left(\int_{\mathbb R^{n+1}_{+}}|(-\Delta)^{\gamma/2}u(x,t)|^{p}t^{p(\gamma-\beta/2)-1}dxdt\right)^{1/p}\approx\|f\|_{\dot{\Lambda}^{p,p}_{\beta/2}(\mathbb R^{n})}\quad\hbox{with}\quad \gamma>\beta/2.
\end{equation}
 The characterizations (\ref{eq-3.1})\ \&\ (\ref{eq-3.2}) indicate that
\begin{equation}\label{equ-3.16}
\left(\int_{\mathbb R^{n+1}_{+}}|\nabla_{(x,t)}u(x,t)|^{p}t^{p-1-p\beta/2}dxdt\right)^{1/p}\approx\|f\|_{\dot{\Lambda}^{p,p}_{\beta/2}(\mathbb R^{n})}, \quad\hbox{when}\quad  \beta\in(0,\min\{2,2s\})
\end{equation}

Equations  (\ref{eq-3.2}) \&\ (\ref{equ-3.15}) \&\   (\ref{equ-3.16}), together with the fractional Sobolev inequality (cf. \cite[Theorem 7.34 ]{AdamsFournier} and \cite[Theorem 6.5]{Di Nezza}):
$$\left(\int_{\mathbb R^{n}}|f(x)|^{np/(n-p\beta/2)}dx\right)^{(n-p\beta/2)/np}\lesssim \|f\|_{\dot{\Lambda}^{p,p}_{\beta/2}(\mathbb R^{n})},\quad  1\leq p<2n/\beta,$$
imply  the following fractional Sobolev trace type inequalities  for $f\in\dot{\Lambda}^{p,p}_{\beta}(\mathbb R^{n})$ with $1\leq p<2n/\beta$,
\begin{equation}\label{eq-3.3}
\|f\|_{L^{\frac{np}{n-p\beta/2}}(\mathbb{R}^n)}
\lesssim \left\{\begin{aligned}
&\left(\int_{\mathbb R^{n+1}_{+}}\left|\nabla_{(x,t)}u(x,t)\right|^{p}t^{p-1-p\beta/2}dxdt\right)^{1/p},&\quad  0<\beta<\min\{2, 2n/p,2s\};
\\
&\left(\int_{\mathbb R^{n+1}_{+}}|\partial_{t}u(x,t)|^{p}t^{p-1-p\beta/2}dxdt\right)^{1/p}, &\quad  0<\beta<\min\{2n/p,2s\};
\\
&\left(\int_{\mathbb R^{n+1}_{+}}|(-\Delta)^{\gamma/2}u(x,t)|^{p}t^{(\gamma-\beta/2)-1}dxdt\right)^{1/p},& \quad  \gamma>\beta/2\ \&\ \beta\in(0, 2n/p).
\end{aligned}\right.
\end{equation}

Moreover,  for $f\in \dot{\Lambda}^{p,p}_{\beta/2}(\mathbb R^{n}),$ there holds the logarithmic type Sobolev inequality (cl. \cite[Corollary 2.6]{{Hajaiej Yu Zhai}}):
$$\exp\left(\frac{\beta}{2n}\int_{\mathbb R^{n}}|f(x)|^{p}\ln|f(x)|^{p}dx\right)\lesssim \|f\|_{\dot{\Lambda}^{p,p}_{\beta/2}(\mathbb R^{n})},\quad \quad \|f\|_{L^{p}(\mathbb R^{n})}=1, \ \beta\in (0,2n/p)\ \&\ p\in(1, 2n/\beta), $$
and the fractional Hardy inequality (cf. \cite[Theorem 1.1]{FS}):
\begin{eqnarray}\label{eq-3.5}
\left(\int_{\mathbb R^{n}}\frac{|f(x)|^{p}}{|x|^{p\beta/2}}dx\right)^{1/p}&\lesssim& \|f\|_{\dot{\Lambda}^{p,p}_{\beta/2}(\mathbb R^{n})},\quad \beta\in(0,2)\ \&\  p\in[1,2n/\beta).
\end{eqnarray}
Similar to (\ref{eq-3.3}),   we can establish the logarithmic and Hardy trace type inequalities for general $p$ with $\beta$ in a similar range.  We omit the details.  Thus, Theorems \ref{them4} \& \ref{them1} \& \ref{them2}   can be generalized to $p\geq 1.$

Note that
$\|f\|_{\dot{H}^{-\beta/2}(\mathbb{R}^n)}\approx\|u(\cdot,\cdot)\|_{L^2(\mathbb{R}^{n+1}_+, \sigma)}$
 is the key to prove Theorem \ref{them3}. For $f\in \dot{\Lambda}^{p,p}_{-\beta/2}(\mathbb{R}^n)$ with  $\beta\in(0,2n)$ and $ p\in(0,\infty), $  Ingmanns proved
  in \cite[Corollary 5.2.3]{Ingmanns}   that \begin{equation}\label{eq3-20}
     \left(\int_{\mathbb R^{n+1}_{+}}|u(x,t)|^{p}t^{p\beta/2-1}dxdt\right)^{1/p}\approx \|f\|_{\dot{\Lambda}^{p,p}_{-\beta/2}(\mathbb{R}^n)}.
 \end{equation}
Thus (\ref{eq3-20}) and  (\ref{affineLp}) imply that  there holds
\begin{eqnarray*}
\|f\|_{\dot{\Lambda}^{p^{\ast},p^{\ast}}_{-\beta/2}(\mathbb{R}^n)}
\lesssim \left(\mathcal{E}_p\left(u,t^{{p^{\ast}\beta}/{2}-1}\right)\right)^{{2n}/{(2n+p^{\ast}\beta)}}\left\|\frac{\partial u}{\partial t}\right\|^{{p^{\ast}\beta}/{(2n+p^{\ast}\beta)}}_{L^p\left(\mathbb{R}^{n+1}_+,t^{{p^{\ast}\beta}/{2}-1}\right)}
\end{eqnarray*}
with $\beta\in (2/p^\ast,2n),$ $1\leq p<n+p^\ast\beta/2$ with $p\ast$ satisfying  $p^\ast\geq \max\{2/\beta,1\}$ and $$1/p=1/p^{\ast}+2/(2n+p^{\ast}\beta).$$ Thus, an inequality similar to (\ref{equ-3.13}) in Theorem \ref{them3} holds for  $ f\in \dot{\Lambda}^{p^{\ast},p^{\ast}}_{-\beta/2}(\mathbb{R}^n).$

\section{Carleson Embeddings of Besov Spaces via the Caffarelli-Silvestre Extension}\label{sec-4}
In the following, we denote the Caffarelli-Silvestre extension of $f$ by $u(\cdot,\cdot).$
In this section, we study the Carleson embeddings via the Caffarelli-Silvestre extension  $u(\cdot,\cdot).$

\subsection{Case (1): $p=q\in(n/(n+\beta),1]$ }

\begin{theorem}
Let $\beta\in (0,n)$ when $p=1$, or  $\beta\in(0,1)$ when $p\in(n/(n+\beta),1).$ Let $q_0\in [p,\infty)$ and $\mu$ be a nonnegative Radon measure on $\mathbb R^{n+1}_{+}$. Then the following statements are equivalent.
\item{\rm (i)} There holds
$$\|u(\cdot,\cdot)\|_{L^{q_0,p}_{\mu}(\mathbb{R}^{n+1}_+)}\lesssim \|f\|_{\dot{\Lambda}^{p,p}_{\beta}(\mathbb{R}^n)}\quad \forall\ f\in C^{\infty}_{0}(\mathbb R^{n}).$$
\item{\rm (ii)} There holds
$$\|u(\cdot,\cdot)\|_{L^{q_0}_{\mu}(\mathbb{R}^{n+1}_+)}\lesssim \|f\|_{\dot{\Lambda}^{p,p}_{\beta}(\mathbb{R}^n)}\quad\forall\ f\in C^{\infty}_{0}(\mathbb R^{n}).$$
\item{\rm (iii)} There holds
$$\|u(\cdot,\cdot)\|_{L^{q_0,\infty}_{\mu}(\mathbb{R}^{n+1}_+)}\lesssim \|f\|_{\dot{\Lambda}^{p,p}_{\beta}(\mathbb{R}^n)}\quad \forall\ f\in C^{\infty}_{0}(\mathbb R^{n}).$$
\item{\rm(iv)} There holds
$$\sup_{t>0}\frac{t^{p/q_0}}{c^\beta_{p,p}(\mu; t)}<\infty.$$
\item{\rm (v)} There holds
$$\sup\left\{\frac{(\mu(T(O)))^{p/q_0}}{C^{p,p}_\beta(O)}: \hbox{open set}\ O\subset \mathbb R^{n}\right\}<\infty.$$
\item{\rm (vi)} There holds
$$\sup_{(x,r)\in\mathbb R^{n+1}_+}\frac{\left(\mu(T(B(x,r)))\right)^{p/q_0}}{r^{n-p\beta}}<\infty.$$
\end{theorem}

\begin{proof} We will prove the equivalence by showing   (i)$\Longrightarrow$(ii)$\Longrightarrow$(iii)$\Longrightarrow$(v)$\Longrightarrow$(i),  (iii)$\Longrightarrow$(iv)$\Longrightarrow$(i), and (v)$\Longleftrightarrow$(vi).

Part 1: (i)$\Longrightarrow$(ii)$\Longrightarrow$(iii)$\Longrightarrow$(v)$\Longrightarrow$(i).
Firstly, note that the Lorentz spaces $L^{p,q}_{\mu}(\mathbb R^{n+1}_{+})$ is increasing with respect to $q$. We can get, as $p\leq q_{0}\leq \infty$,
$$\|u(\cdot,\cdot)\|_{L^{q_{0},\infty}(\mathbb R^{n+1}_{+}, \mu)}\leq \|u(\cdot,\cdot)\|_{L^{q_{0},q_{0}}(\mathbb R^{n+1}_{+}, \mu)}\leq \|u(\cdot,\cdot)\|_{L^{q_{0},p}(\mathbb R^{n+1}_{+}, \mu)}.$$
This gives us  (i) $\Longrightarrow$ (ii) $\Longrightarrow$ (iii).

In the following, we show  (iii)$\Longrightarrow$(v)$\Longrightarrow$(i).
We first show (iii)$\Longrightarrow$(v). Assume that (iii) holds. Let $O\subset \mathbb R^{n}$ be an open set. Let $K$ be a compact subset of $O$. If a nonnegative function $f\in C^{\infty}_{0}(\mathbb R^{n})$ satisfies  $f\geq1$ on $K$, noting that $B(x,t)\subset K$ for $(x,t)\in T(K)$,  we can get for $|x-y|<t$,
$$\frac{c(n,s)t^{s}}{(|x-y|^{2}+t^{2})^{(n+s)/2}}\geq \frac{c(n,s)t^{s}}{(t^{2}+t^{2})^{(n+s)/2}}\geq\frac{c(n,s)}{2^{(n+s)/2}t^{n}}.$$
Hence there exists a constant $c_{n,s}$ depending only on $n$ and $s$ such that
\begin{eqnarray*}
|p^{s}_{t}\ast f(x)|&=&\int_{\mathbb R^{n}}p^{s}_{t}(x-y)f(y)dy\\
&\geq&\int_{K}p^{s}_{t}(x-y)dy\\
&\geq&c(n,s)\int_{B(x,t)}\frac{t^{s}}{(t^{2}+|x-y|^{2})^{(n+s)/2}}dy\\
&\geq &c_{n,s},
\end{eqnarray*}
which indicates that
$$T(K)\subset \Big\{(t,x)\in\mathbb R^{n+1}_{+}:\ |p^{s}_{t}\ast f(x)|>c_{n,s}\Big\}.$$
Hence, by (iii),
\begin{eqnarray*}
\mu(T(K))&\leq&\mu\Big(\Big\{(t,x)\in\mathbb R^{n+1}_{+}:\ |p^{s}_{t}\ast f(x)|>c_{n,s}\Big\}\Big)\\
&\lesssim&\|f\|^{q_0}_{\dot{\Lambda}^{p,p}_{\beta}(\mathbb{R}^n)}.
\end{eqnarray*}
Taking the supermum on $f$ gives
$$\mu(T(K))\leq (C^{p,p}_\beta(K))^{q_0/p}.$$
By the definition of $C^{p,p}_\beta(\cdot)$, we obtain
\begin{eqnarray*}
\mu(T(O))&=&\sup_{K\subset O}\mu(T(K))\\
&\leq&\sup_{K\subset O}(C^{p,p}_\beta(K))^{q_{0}/p}\\
&\leq&(C^{p,p}_\beta(O))^{q_{0}/p}.
\end{eqnarray*}
This proves (v).

Below we consider (v) $\Longrightarrow$ (i).  Given $f\in C^{\infty}_{0}(\mathbb R^{n})$. It follows from (i)-(iv) and (v) of Lemma \ref{lemma 2}  that
$$\mu(E_{\lambda,s})\leq \mu\Big(T\Big(\Big\{x\in\mathbb R^{n}:\ \mathcal{M}f(x)>\lambda/C\Big\}\Big)\Big)$$ with a constant $C$.

If (v) holds for every open set $O\subset \mathbb R^{n}$,
$$\mu(E_{\lambda,s}(f))\lesssim \Big(C^{p,p}_\beta\Big(\Big\{x\in\mathbb R^{n}:\ \mathcal{M}f(x)>\lambda/C\Big\}\Big)\Big)^{q_0/p}.
$$
Then, (i) of Lemma \ref{lemma 24}  implies
\begin{eqnarray*}
\int^{\infty}_{0}\Big(\mu(E_{\lambda,s}(f))\Big)^{p/q_0}d\lambda^p&\lesssim&\int^{\infty}_{0}C^{p,p}_\beta\Big(\Big\{x\in\mathbb R^{n}:\ \mathcal{M}f(x)>\lambda/C\Big\}\Big)d\lambda^p\\
&\lesssim &\|f\|^{p}_{\dot{\Lambda}^{p,p}_{\beta}(\mathbb{R}^n)},
\end{eqnarray*}
i.e., (i) holds.

Part 2: (iii)$\Longrightarrow$(iv)$\Longrightarrow$(i).
If (iii) holds, for any bounded open set $O\subset\mathbb R^{n}$, let $f\in C^{\infty}_0(\mathbb{R}^n)$ with $f\geq 1$ on $O$. For $(x,t)\in T(O)$, then $p^{s}_{t}\ast |f|(x)>\eta_{n,s}$. Then
\begin{eqnarray*}
\left(\mu(T(O))\right)^{1/q_0}&\leq&\Big(\mu\Big(\Big\{(x,t)\in\mathbb R^{n+1}_{+}:\ p^{s}_{t}\ast f(x)>\eta_{n,s}\Big\}\Big)\Big)^{1/q_0}\\
&\leq&K_{p}(\mu)\|f\|_{\dot{\Lambda}^{p,p}_{\beta}(\mathbb{R}^n)}.
\end{eqnarray*}
For any $0<t<\mu(T(O))$, we have $t\leq (K_{p}(\mu)\|f\|_{\dot{\Lambda}^{p,p}_{\beta}(\mathbb{R}^n)})^{q_0}$. Taking the infimum on $\|f\|_{\dot{\Lambda}^{p,p}_{\beta}(\mathbb{R}^n)}$ reaches
$t\leq(K_{p}(\mu))^{q_0}(C^{p,p}_\beta(O))^{q_0/p}$. This gives
$$\sup_{t>0}\frac{t^{p/q_0}}{c^\beta_{p,p}(\mu; t)}\leq (K_{p}(\mu))^{p}.$$
On the other hand, if (iv) holds, then we can get
\begin{eqnarray*}
&&\int^{\infty}_{0}\mu(E_{\lambda,s}(f)\cap T(B(0,k)))^{p/q_0}d\lambda^{p}\\
&&\leq\int^{\infty}_{0}\left\{\frac{(\mu(E_{\lambda,s}(f)\cap T(B(0,k))))^{p/q_0}}{c^\beta_{p,p}(\mu; \mu(E_{\lambda,s}(f)\cap T(B(0,k))))}\right\}\times C^{p,p}_\beta(O_{\lambda,s}(f)\cap B(0,k))d\lambda^{p}\\
&&\leq \left(\sup_{t>0}\frac{t^{p/q_0}}{c^\beta_{p,p}(\mu;t)}\right)\int^{\infty}_{0}C^{p,p}_\beta\Big(\Big\{x\in\mathbb R^{n}:\ \theta_{n,s} \mathcal{M}f(x)>\lambda\Big\}\cap B(0,k)\Big)d\lambda^{p}\\
&&\leq \left(\sup_{t>0}\frac{t^{p/q_0}}{c^\beta_{p,p}(\mu;t)}\right)\|f\|^{p}_{\dot{\Lambda}^{p,p}_{\beta}(\mathbb{R}^n)}.
\end{eqnarray*}
Letting $k\rightarrow \infty$ gives us (i).

Part 3: (v)$\Longleftrightarrow$(vi). We first assume that (v) holds. Then let $O=B(x,r)$. By (v),
$$(\mu(T(B(x,r))))^{p/q_0}\lesssim C^{p,p}_\beta (B(x,r)).$$

For $d\in(0,\infty)$, the classical Hausdorff capacity $H^{(\infty)}_{d}$ is determined via
$$H^{(\infty)}_{d}(K):=\inf\Big\{\sum^{\infty}_{i=0}m_{i}2^{-di}\Big\},$$
where the infimum is taken over all countable coverings of $K\subset \mathbb R^{n}$ by balls whose radii $r_{j}<\infty$, while $m_{i}$ is the number of balls from this covering whose radii $r_{j}$ belong to the interval $(2^{-i-1}, 2^{-i}], i=0,1,2,\ldots$.
According to \cite[Theorem 3.6, Theorem 3.8]{Adams1} and \cite[page 32]{Adams 1998}, we have
$C^{p,p}_\beta(\cdot)\approx H^{(\infty)}_{n-\beta p}(\cdot).$ Since $C^{p,p}_\beta(B(x,r))\approx H^{(\infty)}_{n-p\beta}(B(x,r))\approx r^{n-p\beta}$,
\begin{eqnarray*}
\sup_{(x,r)\in \mathbb{R}^{n+1}_+}\frac{(\mu(T(B(x,r))))^{p/q_0}}{r^{n-p\beta}}&\approx& \sup_{(x,r)\in \mathbb{R}^{n+1}_+}\frac{(\mu(T(B(x,r))))^{p/q_0}}{C^{p,p}_\beta(B(x,r))}\\
&\lesssim& \sup\left\{\frac{(\mu(T(O)))^{p/q_0}}{C^{p,p}_\beta(O)}:\ \hbox{open set}\ O\subset \mathbb R^{n}\right\}<\infty.
\end{eqnarray*}
This gives (vi).

Conversely, we assume that (vi) holds.  Then we have
$$T_{q_0,n,\beta,p}(\mu):=\sup_{(x,r)\in \mathbb{R}^{n+1}_+}\frac{\mu\left(T(B(x,r))\right)}{r^{q_0(n-p\beta)/p}}<\infty.$$
Given an open set $O\subset\mathbb R^{n}$ and let $\{B(x_{j}, r_{j})\}_{j\in\mathbb Z_{+}}$ be a sequence of covering balls of $O$, i.e., $O\subseteq \cup_{j}B(x_{j}, r_{j})$.   Below we only need to verify that (vi) holds for bounded sets. In fact, if (v) holds for bounded sets, for any open set $O\subset\mathbb{R}^n,$ let $E_{0}=B(0,1)$ and $E_{k}=B(0, 2^{k})\setminus B(0, 2^{k-1})$, $k\geq 1$. Then $\mathbb R^{n}=\cup^{\infty}_{k=0}E_{k}$ and
$$O=O\cap\Big(\cup^{\infty}_{k=0}E_{k}\Big)=\cup^{\infty}_{k=1}\Big(O\cap E_{k}\Big).$$
So, we can get
\begin{eqnarray*}
(\mu(T(O)))^{p/q_0}
&\leq&\left(\mu(\cup^{\infty}_{k=1}T( (O\cap E_{k} )))\right)^{p/q_0}\\
&\leq&\left(\sum^{\infty}_{k=0}\mu(T(O\cap E_{k}))\right)^{p/q_0}\\
&\lesssim&\left(\sum^{\infty}_{k=0}C^{p,p}_\beta(O\cap E_{k})^{q_0/p}\right)^{p/q_0}\\
&\lesssim&\sum^{\infty}_{k=0}C^{p,p}_\beta(O\cap E_{k})\\
&\lesssim&C^{p,p}_\beta(\cup^{\infty}_{k=0}(O\cap E_{k}))\\
&\lesssim&C^{p,p}_\beta(O).
\end{eqnarray*}

In the following, we assume that $O$ is a bounded open subset of $\mathbb{R}^n.$
Then, there exists a dyadic cube sequence $\{I_{j}^{(1)}\}$ in $\mathbb R^{n}$ such that
$O\subseteq\cup_{j}I_{j}^{(1)}$ and $|B(x_{j}, r_{j})|\approx |I_{j}^{(1)}|$. On the other hand,  according to Dafni-Xiao \cite[Lemma 4.1]{DafniXiao}, there exists a sequence of dyadic cubes $\{I^{(2)}_{j}\}$ and denote by $\{I^{(3)}_{j}\}$ its $5\sqrt{n}$ expansion such that
$$\left\{\begin{aligned}
&\cup_{j}I^{(2)}_{j}=\cup_{j}I^{(1)}_{j};\\
&\sum_{j}\left|I_{j}^{(2)}\right|^{1-p\beta/n}\leq \sum_{j}\left|I^{(1)}_{j}\right|^{1-p\beta/n};\\
&T(O)\subseteq \cup_{j}T(I^{(3)}_{j}).
\end{aligned}
\right.$$
Since $|B(x_{j}, r_{j})|\approx|I^{(1)}_{j}|$, $r_{j}\approx\ell(I^{(1)}_{j})\approx \ell(I^{(3)}_{j})$. Let $x_{I_{j}}$ be the center of $I^{(1)}_{j}$. Then $I^{(3)}_{j}\subseteq B(x_{I_{j}}, \ell(I_{j}^{(3)}))$. We can get
\begin{eqnarray}\label{them1Inequ 4}
\mu(T(O))&\leq&\sum_{j}\mu(T(I^{(3)}_{j}))\\
&\leq&\sum_{j}\mu\left(T(B(x_{I_{j}}, \ell(I_{j}^{(3)})))\right)\nonumber\\
&\leq&\sum_{j}\left(\ell(I^{(3)}_{j})\right)^{q_0(p\beta-n)/p}\mu\left(T(B(x_{I_{j}}, \ell(I_{j}^{(3)})))\right)\left(\ell(I^{(3)}_{j})\right)^{q_0(n-p\beta)/p}\nonumber\\
&\lesssim&T_{q_0,n,\beta,p}(\mu)\sum_{j}\left(\ell\left(I^{(3)}_{j}\right)\right)^{q_0(n-p\beta)/p}\nonumber\\
&\lesssim&T_{q_0,n,\beta,p}(\mu)\sum_{j}\Big|I^{(3)}_{j}\Big|^{q_0(1-p\beta/n)/p}.\nonumber
\end{eqnarray}
On the other hand, we have
\begin{eqnarray*}
\sum_{j}\left|I^{(3)}_{j}\right|^{q_0(1-p\beta/n)/p}
\lesssim\sum_{j}\left|I^{(2)}_{j}\right|^{q_0(1-p\beta/n)/p}
\lesssim\Big(\sum_{j}\left|I^{(1)}_{j}\right|^{1-p\beta/n}\Big)^{q_0/p}
\lesssim\Big(\sum_{j}\left|B(x_{j}, r_{j})\right|^{1-p\beta/n}\Big)^{q_0/p},
\end{eqnarray*}
which implies
\begin{eqnarray}\label{them1Inequ 5}
\sum_{j}\left|I^{(3)}_{j}\right|^{q_0(1-p\beta/n)/p}
\lesssim\left(H^{(\infty)}_{n-p\beta}(O)\right)^{q_0/p}
\lesssim\left(C^{p,p}_\beta(O)\right)^{q_0/p}
\end{eqnarray}
according to the definition of $H^{(\infty)}_{n-p\beta}$ and the equivalence of $C^{p,p}_\beta(\cdot)\approx H^{(\infty)}_{n-\beta p}(\cdot).$
Thus, by (\ref{them1Inequ 4}) and (\ref{them1Inequ 5}),  we get  $\mu(T(O))\lesssim T_{q_0,n,\beta,p}(\mu)\left(C^{p,p}_\beta(O)\right)^{q_0/p}$ which implies (v).
\end{proof}

\subsection{Case (2): $(p,q)\in (1,n/\beta)\times (1,\infty)$}

Denote $p\lor q=\max\{p,q\}$ and $p\land q=\min\{p,q\}.$

\begin{theorem}\label{th-5.1}
Let $\mu$ be  a nonnegative Borel measure on $\mathbb R^{n+1}_{+}$. If  $\beta\in (0,n)$  and $(p,q)\in(1,{n}/{\beta})\times(1,\infty),$ then
\begin{equation}\label{eq-5.1}
\|u(\cdot,\cdot)\|_{L^{p\land q,p\lor q}_{\mu}(\mathbb{R}^{n+1}_+)}
\lesssim \|f\|_{\dot{\Lambda}^{p,q}_{\beta}} \quad \forall\ f\in C^{\infty}_0(\mathbb{R}^n) \end{equation}
if and only if
\begin{equation}\label{eq-5.2}
\mu(T(O))\lesssim (C^{p,q}_\beta(O))^{{(p\land q)}/{p}}\quad \forall  \
 \text{ open set } O\subset \mathbb{R}^n.
\end{equation}

\end{theorem}
\begin{proof}
We first assume that (\ref{eq-5.2}) holds. Recall that
$$O_{\lambda,s}(f)=\Big\{x\in\mathbb R^{n}:\ \sup_{|y-x|<t}|u(y,t)|>\lambda\Big\},\quad E_{\lambda,s}(f)=\Big\{(x,t)\in\mathbb R^{n+1}_{+}:\ |u(x,t)|>\lambda\Big\}.$$
It follows from  (i) of Lemma \ref{lemma 2} that $$E_{\lambda,s}(f)\subset T(O_{\lambda,s}),\quad \mu(E_{\lambda,s}(f))\leq \mu(T(O_{\lambda,s}(f))).$$ So,   (\ref{eq-5.2}) implies
$$\mu((E_{\lambda,s}(f))\leq \mu(T(O_{\lambda,s}(f)))\\
\lesssim (C^{p,q}_\beta(O_{\lambda,s}(f))^{{(p\land q)}/{p}},
$$
which gives
\begin{eqnarray*}
\int^{\infty}_{0}(\mu(E_{\lambda,s}(f)))^{{(p\lor q)}/{(p\land q)}}d\lambda^{p\lor q} &\leq&\int^{\infty}_{0}(\mu(T(O_{\lambda,s}(f))))^{{(p\lor q)}/{(p\land q)}}d\lambda^{p\lor q}\\
&\lesssim&\int^{\infty}_{0}C^{p,q}_\beta(O_{\lambda,s}(f))^{1\lor ({q}/{p})}d\lambda^{p\lor q}.
\end{eqnarray*}
It follows from (iv) and (v) of Lemma \ref{lemma 2}  that
$O_{\lambda,s}\subset \Big\{x\in\mathbb R^{n}:\ \mathcal{M}f(x)>\lambda/C\Big\}$ with a constant $C$. So, we get
$$C^{p,q}_\beta(O_{\lambda,s})\leq C^{p,q}_\beta\Big(\Big\{x\in\mathbb R^{n}:\ \mathcal{M}f(x)>\lambda/C\Big\}\Big),$$
which implies that
\begin{eqnarray*}
\int^{\infty}_{0}(\mu(E_{\lambda,s}(f)))^{{(p\lor q)}/{(p\land q)}}d\lambda^{p\lor q}&\lesssim& \int^{\infty}_{0}C^{p,q}_\beta(\{x\in\mathbb R^{n}:\ \mathcal{M}f(x)>\lambda/C\})^{1\lor ({q}/{p})}d\lambda^{p\lor q}\\
&\lesssim& \|f\|^{p\lor q}_{\dot{\Lambda}^{p,q}_{\beta}(\mathbb{R}^n)}
\end{eqnarray*}
according to (i) of Lemma \ref{lemma 24}. Hence (\ref{eq-5.1}) holds.

\end{proof}

\subsection{Case (3):  $(p,q)\in (1,n/\beta)\times \{\infty\}$}
\begin{theorem}\label{them sec5}
Let  $\beta\in(0,1)$ and $\mu$ be a non-negative outer measure on $\mathbb{R}^{n+1}_+.$
 If $(\beta, p, q_0)\in (0,1)\times (1,{n}/{\beta})\times (0,\infty),$ then
\begin{equation}\label{them 4.3-1}
\|u(\cdot,\cdot)\|_{L^{q_0,\infty}_{\mu}(\mathbb{R}^{n+1}_+)}\lesssim \|f\|_{\dot{\Lambda}^{p,\infty}_{\beta}(\mathbb{R}^n)}\quad \forall\ f\in C^{\infty}_0(\mathbb{R}^n)
\end{equation}
if and only if
\begin{equation} \label{them 4.3-2}
\mu(T(O))\lesssim (C^{p,\infty}_\beta(O))^{q_0/p}\quad \forall \  \hbox{open set}\   O\subset \mathbb{R}^n.
\end{equation}

\end{theorem}
\begin{proof}  Let $(\beta, p, q_0)\in (0,1)\times (1, {n}/{\beta})\times (0,\infty).$ Assume
(\ref{them 4.3-1}) holds.
Then $$K:=\sup_{0<\|f\|_{\dot{\Lambda}^{p,\infty}_{\beta}(\mathbb{R}^n)}<\infty, f\in C^{\infty}_0(\mathbb{R}^n)}\frac{\sup_{\lambda\in(0,\infty)}
\left(\lambda^{q_0}\mu\left(E_{\lambda,s}(f)\right)\right)^{1/q_0}}{\|f\|_{\dot{\Lambda}^{p,\infty}_{\beta}(\mathbb{R}^n)}}<\infty.$$
For  any open set $O\subset \mathbb{R}^n,$ we can take  $f\in C^{\infty}_0(\mathbb{R}^n)$ such that $O\subset\hbox{Int}(\{x\in \mathbb{R}^n: f(x)\geq 1\})$.  By (iii) of Lemma \ref{lemma 2},
$$(\mu(T(O)))^{1/q_0}\lesssim K\|f\|_{\dot{\Lambda}^{p,\infty}_{\beta}(\mathbb{R}^n)}.$$
Then the definition of $C^{p,\infty}_\beta(\cdot)$ implies $$\mu(T(O))\lesssim (C^{p,\infty}_\beta(O))^{q_0/p}. $$ Thus, (\ref{them 4.3-2}) holds.

On the other hand, assume that (\ref{them 4.3-2}) holds. Then for   any $f\in C^{\infty}_0(\mathbb{R}^n),$  (i)-(ii) of Lemma \ref{lemma 2} imply
\begin{eqnarray*}
\mu(E_{\lambda,s}(f))&\lesssim&\mu(T(\{x\in \mathbb{R}^{n}: \mathcal{M}f(x)>\theta_{n,s}\lambda\}))\\
&\lesssim&(C^{p,\infty}_\beta(T(\{x\in \mathbb{R}^{n}: \mathcal{M}f(x)>\theta_{n,s}\lambda\})))^{q_0/p}.
\end{eqnarray*}
So, (ii) of Lemma \ref{lemma 24} indicates that
$$(\lambda^{q_0}\mu(E_{\lambda,s}(f)))^{1/q_0}\lesssim(\lambda^pC^{p,\infty}_\beta(T(\{x\in \mathbb{R}^{n}: \mathcal{M}f(x)>\theta_{n,s}\lambda\})))^{1/p}\lesssim\|f\|_{\dot{\Lambda}^{p,\infty}_{\beta}(\mathbb{R}^n)}.
$$
We have
$$\sup_{\lambda>0}
\left(\lambda^{q_0}\mu\left(E_{\lambda,s}(f)\right)\right)^{1/q_0}\lesssim \|f\|_{\dot{\Lambda}^{p,\infty}_{\beta}(\mathbb{R}^n)},$$
which gives (\ref{them 4.3-1}).
\end{proof}

\section*{Acknowledgement}
The authors are very grateful to the anonymous reviewers for their  numerous useful suggestions which lead to a substantial improvement of this manuscript.

This work was in part supported by the  National Natural Science Foundation of
China (\# 11871293 \ \&\ \# 12071272) \&  Natural Science
Foundation of Shandong Province (No.  ZR2020MA004).

\end{document}